\documentclass{article}
\usepackage[left=3.5cm, top=2cm, right=2.5cm, bottom=20mm, nohead, nofoot]{geometry}
\title{Cyclic prime numbers}
\author{Konstantin Kutsenko}
\date{May 2021}

\usepackage{aeguill}
\usepackage{graphicx}
\usepackage{amsmath}
\usepackage{wrapfig}
\usepackage{float}

\usepackage[sans]{dsfont}

\graphicspath{{pictures/}}
\DeclareGraphicsExtensions{.pdf,.png,.jpg}

\begin{document}

\maketitle

\begin{abstract}
\noindent This work is meant to demonstrate new class of prime numbers - cyclic prime numbers, that can be derived from any prime number at certain numeric systems. Cyclic prime numbers are also related to the cyclic numbers and full reptend prime numbers. Cyclic prime numbers, derived from the same prime number, could be grouped. Conditions and properties of those groups between different numeric systems are also the subjects of this paper.
\end{abstract}

\paragraph{1. Introduction} ~\\

\noindent First publication related to cyclic numbers and full reptend prime numbers appears in the book
[1].

\bigskip

\noindent This book was published 200 years before this work, authors of this book show that periodic fraction derived from 1/P can be presented as converging geometric series. Authors don't use the term geometric series, but they show that 1/7 is the sum of decreasing rational numbers.

\bigskip

\noindent In the book [2] there is formulation of a condition under which full reptend prime numbers appear.

\bigskip

\noindent In book [3] there is a mentioning of cyclic numbers and their relation to repunits.

\bigskip

\noindent In the book [4] there is a mentioning of divisions during the formation of periodic fractions, which are subsequently used in the formula of geometric series.

\bigskip

\noindent In the above mentioned books, and throughout the Internet I wasn't able to find all the regularities that I came up with, so I wish to share them in this article.

\paragraph{2. Cyclic numbers}  ~\\

\noindent A cyclic number is an integer in which cyclic permutations of the digits are successive integer multiples of the number.

\noindent The most famous cyclic number is 142857. It is popularized with pseudoscientific theory of an 'enneagram', however there are not so many scientific results can be found about it.
Let's take a look at its cyclic property:

\begin{center}
142857 * 2 = 285714

142857 * 3 = 428571

142857 * 4 = 571428

142857 * 5 = 714285

142857 * 6 = 857142
\end{center}

\noindent As we can see, when the original number 142857 is multiplied by the numbers from 2 to 6, we get
cyclic permutations of the number 142857.

\bigskip

\noindent Those numbers have other properties, for example regularities that can be observed when multiplying by numbers greater than 7, but they are not a subject of this article.

\paragraph{3. Cyclic prime numbers} ~\\

\noindent A cyclic prime is a prime number formed from a sequence of digits in a cyclic number that repeats for more than one cycle.

\bigskip

\noindent The first cyclic prime formed from 142857 is 1428571, which is a prime number.
Such a number can be written by the first digit of the initial cyclic number and the total number of digits. For example, for 1428571, the first digit is 1 and the total number of digits is 7.

\bigskip

\noindent Here are all the primes formed from the cyclic number 142857, and not exceeding titanic primes (up to 10 thousand digits). The first numbers are written in full, the longer ones
described by the first digit and the number of digits.

\bigskip

\noindent The first 7 cyclic primes, formed from 142857:
1428571, 
71428571, 
7142857142857, 
571428571428571, \linebreak
1428571428571428571428571,
28571428571428571428571428571,
7142857142857142857142857142857.

\noindent Amount of digits in those numbers: 7, 8, 13, 15, 25, 29, 31.

\bigskip

\noindent Rest of the numbers are written by first digits and total amount of digits:

\bigskip

\begin{tabular}{ l l}
First digit & Amount of digits  \\
2 & 34 \\
4 & 41 \\
7 & 104 \\
5 & 273  \\
2 & 304 \\
1 & 355 \\
7 & 440 \\
7 & 571 \\
1 & 823  \\
7 & 2215 \\
5 & 2523 \\
4 & 4379 \\
2 & 4510 \\
4 & 7553 \\
4 & 7679 \\
7 & 9536 \\
\end{tabular}

\bigskip

\noindent In total there are 23 cyclic primes, that are not greater than $10^{1000}$.

\paragraph{4. Properties of prime numbers dependent on numeric system. Full reptend prime numbers} ~\\

\noindent In order to consider how do cyclic primes appear, we also need to consider how cyclic numbers appear.

\bigskip

\noindent There is a class of prime numbers called full reptend prime.
This class of prime numbers depends on the numeric system, so any prime number is a full reptend at certain numeric system.

\bigskip

\noindent P is a prime number. If the periodic fraction formed during the calculation of the rational number 1 / P, has a period equal to P-1 in some numeric system N, then we can say that the prime number P in the numeric system N is full reptend.

\bigskip

\noindent If the number P is a full reptend in some numeric system N, then all P-1 of its digits form
a cyclic number.
Consider a prime number P = 7 in decimal notation.
The number formed from 1 / 7 = $0.\overline{142857}$. The period is 6, which is equal to P-1.

\bigskip

\noindent Consider the remaining fractions from the set 1 / P .. P-1 / P:

2 / 7 = $0.\overline{285714}$

3 / 7 = $0.\overline{428571}$

4 / 7 = $0.\overline{571428}$

5 / 7 = $0.\overline{714285}$

6 / 7 = $0.\overline{857142}$

\bigskip

\noindent We observe that those fractions have the property of a cyclic permutation.

\bigskip

\noindent One of the visualizations will be presented below. In the table below, each row is
a prime number.

\noindent Each column represents a numeric system. The value in a cell is the length of the period of the rational 1 / P. Highlighted in the green
cells that are full reptend.

\bigskip

\noindent First, let's take a look at the individual parts of this table.
For each prime number P, there is a sequence of possible period lengths of the fraction 1 / P.

\bigskip

\begin{figure}[h]
\noindent\includegraphics{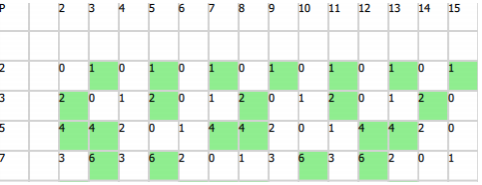}
\caption{1/P in first 14 numeric systems for P = 2, 3, 5, 7}
\label{fig:image}
\end{figure}

\begin{figure}[h]
\noindent\includegraphics{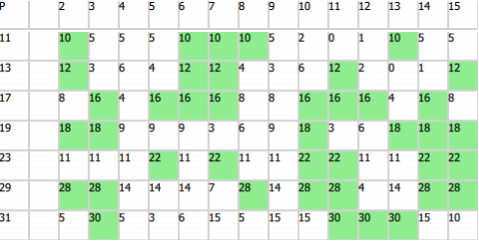}
\caption{1/P in first 14 numeric systems for P = 11, 13, 17, 19, 23, 29, 31}
\label{fig:image}
\end{figure}

\bigskip

\noindent For P = 2, the numeric system cycle length is 2.
\noindent For P = 3, the numeric system cycle length is 3, etc.

\bigskip

\noindent To calculate the length of a period in a certain numeric system, which is represented as base - we need to solve the  equation:

\begin{equation}
  {base}^{period}\bmod P=1
\end{equation}

\bigskip

\noindent If the base is coprime to the P then Fermat's little theorem says that Fermat quotient is an integer. If the base is also a generator of multiplucative group of integers modulo p, then Fermat quotient will be cyclic number, and p will be a full reptend prime.

\bigskip

\noindent Let's take a look at the figure 2 for different P in different numeric systems.

\bigskip

\noindent In the decimal system, we see that the first prime numbers, which are full reptends, are: 7, 17, 19, 23, 29.

\bigskip

\noindent The numbers 2 and 5 do not have a period here, since the base of the numeric system is divisible by both of these prime numbers without remainder.

\bigskip

\noindent In the case of P = 3, we get a periodic fraction with a unit length of the period in the decimal system: 1/3 = 0,(3).

\bigskip

\noindent In the case of P = 11, we get a periodic fraction with a period length of 2 in decimal: 1/11 = 0, (09).

\bigskip

\noindent With P = 13 we get a special case, the period length is 6, which is not equal to P-1.
However, the period length is (P-1) / 2. When this proportion is respected, two sets of cyclic numbers are formed.
Such P is called the 2nd reptend level prime.

\bigskip

\noindent Here's the example of 2nd reptend level prime:

\begin{center}
\Large$\frac{1}{13} = 0.\overline{076923}$ \vspace{0.1cm}

\Large$\frac{2}{13} = 0.\overline{153846}$ \vspace{0.1cm}
\end{center}

\noindent All other fractions from P = 13 up to P-1 / P will have the same digits as 1/13 or 2/13, but with cyclic permutations.

\begin{center}

\Large$\frac{3}{13} = \normalsize 0.\overline{230769}$  \vspace{0.1cm}

\Large$\frac{4}{13} = \normalsize 0.\overline{307692}$ \vspace{0.1cm}

\Large$\frac{5}{13} = \normalsize 0.\overline{384615}$ \vspace{0.1cm}

\Large$\frac{6}{13} = \normalsize 0.\overline{461538}$ \vspace{0.1cm}

\Large$\frac{7}{13} = \normalsize 0.\overline{538461}$ \vspace{0.1cm}

\Large$\frac{8}{13} = \normalsize 0.\overline{615384}$ \vspace{0.1cm}

\Large$\frac{9}{13} = \normalsize 0.\overline{692307}$ \vspace{0.1cm}

\Large$\frac{10}{13} = \normalsize 0.\overline{769230}$ \vspace{0.1cm}

\Large$\frac{11}{13} = \normalsize 0.\overline{846153}$ \vspace{0.1cm}

\Large$\frac{12}{13} = \normalsize 0.\overline{923076}$ \vspace{0.1cm}

\end{center}

\noindent N-reptend prime numbers also form cyclic primes: primes can be formed from each of the sequences of cyclic numbers.

\bigskip

\noindent There are prime numbers for the first cycle: 769230769, 769230769230769230769. But also there is a prime number from the second cycle: 1538461.

\bigskip

\noindent Concluding this topic it is interesting to note that the numeric systems in which we meet full reptend primes are also unusual.

\bigskip

\noindent For P = 7, the first 2 numeric systems, in which it is a full reptend, are systems with N 3 and 5 - twin primes.
Further, these positions are repeated every 7 numeric systems, and always remain in same distance. When the sum of the bases of the numeric systems is evenly divisible by 12, we meet
twin primes again. For example, 17 and 19, 59 and 61.

\paragraph{5. Representation of a periodic fraction in the form of a converging geometric series} ~\\

\noindent Each of the fractions formed by the full reptend or n-th repntend level can be decomposed into converging geometric series. For every prime P and a given numeric system N, there is an infinite number of such geometric series.

\bigskip

\noindent Formula for writing the sum of a geometric progression for 1 / P is:

\begin{equation}
  \sum\limits_{n=0}^\infty\frac{s*r^n}{base^{length(n+1)}} = \frac{1}{P}
\end{equation}

\noindent Where s is an integer derived from the fraction 1 / P.
Since full reptend prime form an infinite periodic fraction, we can get the set of s from its numbers with the amount of digits from 1 to infinity:

\begin{equation}
	s=[\frac{1}{P}*base^{length}]
\end{equation}

\noindent The length parameter denotes the number of symbols that will be used in the number s, this parameter can be varied from one to infinity. For each new length parameter, we will get a new geometric progression.

\bigskip

\noindent The number r is also an integer and is one of the remainders formed by calculating 1 / P. As a fraction, 1 / P will have a period P-1 if the prime number is full reptend in the numeric system, in the same way the number of unique remainders will be equal to P-1.

\begin{wrapfigure}[15]{l}{0.35\linewidth} 
\includegraphics{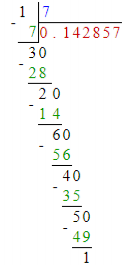}
\end{wrapfigure}

\bigskip

\noindent In order to understand how division remainders are formed we need to use column division.

\bigskip

\noindent As a result, we got the remainders: [3, 2, 6, 4, 5, 1]. It is these values that will take part in the formula. Let's derive mathematical formula.

\bigskip

\noindent The first remainder of the division will be base mod P.

\noindent Each subsequent remainder will depend on the previous one, the function can be written recursively:

\bigskip

\begin{equation}
 \begin{cases}
   r_0 = 1 \\
   r_n = r_{n-1} * (base\mod P) \\
 \end{cases}
\end{equation}

\bigskip

\noindent Let's translate the recursive formula into a closed form:

\begin{equation}
r_{length}=base^{length}\mod P
\end{equation}

\bigskip

\noindent Let's write down the general formula of the geometric series using only the following parameters: P
- prime number; base - investigated numeric system; length is a parameter that determines the number of the geometric series for a given prime number and a given numeric system.

\begin{equation}
\frac{1}{P} = \sum\limits_{n=0}^\infty\frac{[\frac{1}{P}*base^{length}]*(base^{length}\mod P)^n}{base^{length(n+1)}}
\end{equation}

\bigskip

\noindent Here are the formulas for P = 7 using different s, starting with the shortest ones:

\begin{equation}
\frac{1}{7} = \sum\limits_{n=0}^\infty\frac{1*3^n}{10^{n+1}}
\end{equation}

\noindent In this formula, s = 1, this is the first digit from the fraction 0, (142857), i.e. parameter length = 1.

\noindent The remainder is r = 3, this is the very first remainder, which corresponds to the parameter length = 1.

\bigskip

\begin{equation}
\frac{1}{7} = 0.1 + 0.03 + 0.009 + 0.0027 + 0.00081 + .. 
\end{equation}

\noindent Each next term in the progression is obtained by multiplying by 3 and dividing by 10.

\begin{equation}
\frac{1}{7} = \sum\limits_{n=0}^\infty\frac{14*2^n}{10^{2(n+1)}}
\end{equation}

\begin{equation}
\frac{1}{7} = 0.14 + 0.0028 + 0.000056 + 0.00000112 + ..
\end{equation}

\noindent In this formula, each next term of the progression is obtained by multiplication by 2 and
division by 100.
Here s = 14, these are the first two digits from the fraction $0.\overline{142857}$, i.e. parameter length = 2.
In this case, the remainder is r = 2, this is the second remainder, which corresponds to the parameter length = 2.

\begin{figure}[H]
\begin{center}
\noindent\includegraphics{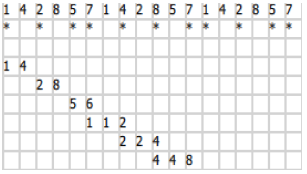}
\caption{Geometric progression in 1/7}
\label{fig:image}
\end{center}
\end{figure}

\noindent There are interesting patterns that can be explored by observing this geometric progression presented in figure 3.

\bigskip

\noindent Places of occurrence of a new element of the progression, where the first digit is not equal to zero (asterisks in the figure above), when investigating and trying to describe the periodicity, showed a rather interesting result - the are proportions that correspond to the major music modes like Ionian, Locrian, Lydian etc. And seams like there is no period, so the sequence is infinite and unique, like digits sequence of $\pi$. This is a very interesting topic, but it has nothing to do with cyclic primes and it deserves a separate article.

\bigskip

\noindent Further formulas, by analogy with the previous ones, can be obtained simply by sequentially increasing the value of length by 1.

\begin{equation}
\frac{1}{7} = \sum\limits_{n=0}^\infty\frac{142*6^n}{10^{3(n+1)}}
\end{equation}

\begin{equation}
\frac{1}{7} = \sum\limits_{n=0}^\infty\frac{1428*4^n}{10^{4(n+1)}}
\end{equation}

\begin{equation}
\frac{1}{7} = \sum\limits_{n=0}^\infty\frac{14285*5^n}{10^{5(n+1)}}
\end{equation}

\begin{equation}
\frac{1}{7} = \sum\limits_{n=0}^\infty\frac{142857*1^n}{10^{6(n+1)}}
\end{equation}

\noindent Finally, we get a sequence in which s is a prime number:

\begin{equation}
\frac{1}{7} = \sum\limits_{n=0}^\infty\frac{1428571*3^n}{10^{7(n+1)}}
\end{equation}

\noindent At this point, consideration of further geometric series can be postponed, but it is important to 
notice that there can be an infinite set of them, since the numbers s for each number P in
some numeric system N is an infinite set:

\begin{equation}
N \in \mathds{N}
\end{equation}

\bigskip

\noindent Let's take a look at geometric series for P = 17:

\begin{equation}
\frac{1}{17} = \sum\limits_{n=0}^\infty\frac{5*15^n}{10^{2(n+1)}}
\end{equation}

\begin{equation}
\frac{1}{17} = \sum\limits_{n=0}^\infty\frac{58*14^n}{10^{3(n+1)}}
\end{equation}

\begin{equation}
\frac{1}{17} = \sum\limits_{n=0}^\infty\frac{588*4^n}{10^{4(n+1)}}
\end{equation}

\noindent It is interesting to consider the decomposition of the number 89 into geometric series.
1/89 = 0.0112359 .. - you can see how the Fibonacci numbers are observed in the first digits of the fraction:

\begin{equation}
\frac{1}{89} = \sum\limits_{n=0}^\infty\frac{1*11^n}{10^{2(n+1)}}
\end{equation}

\begin{equation}
\frac{1}{89} = \sum\limits_{n=0}^\infty\frac{Fibonacci(n)}{10^{n+1}}
\end{equation}

\noindent Interestingly, a similar phenomenon can be found in another prime number - 109:

\bigskip

\begin{equation}
\frac{1}{109} = \sum\limits_{n=0}^\infty\frac{9*17^n}{10^{3(n+1)}}
\end{equation}

\begin{equation}
\frac{1}{109} = \sum\limits_{n=0}^\infty\frac{Fibonacci(n)*(2(n\mod2)-1))}{10^{n+1}}
\end{equation}

\bigskip

\noindent The formula for this number differs from the formula 1/89 only by a factor: 2 * (n mod 2) - 1. This multiplier leads to the fact that it is not summation that occurs, but an alternation of summation and subtraction different elements of the progression.

\bigskip

\noindent Perhaps there are other expansions of the periodic fraction, in addition to the converging
geometric series and summation of Fibonacci numbers.

\paragraph{6. Cyclic and sub-cyclic primes formation} ~\\

\noindent When we consider different values of s from the geometric series formula, we
note that there will be prime numbers among them.

\bigskip

\noindent For example, for a prime number P = 7 forming a cyclic number 142857, the first prime the number will be 1428571. Further, a large number of similar numbers will appears. In order to consider all possible variate of them, it is necessary to consider not only geometric series of 1 / P, but also all other progressions from the set 1 / P .. P-1 / P. Otherwise we might skip, for example, the prime number 71428571.

\bigskip

\noindent Presumably, there is an infinite set of such numbers in each numeric system. However, I could not find the pattern of their occurrence. Perhaps, such numbers are of a limited quantity, and if so, in my opinion it would be even more interesting, but I assume that the set is infinite.

\bigskip

\noindent In fact, if we look closely at the parameter s, we will see prime numbers in it,
which will be less than a cyclic number, but will contain a sequence of digits
from the cyclic number. Such numbers can be called simple subcyclic numbers.

\bigskip

\noindent Let's consider them using the example P = 7. If we simultaneously consider not only 1 / P, but also
all other fractions up to P-1 / P, then we will see that among the parameters s there is a set of 
other prime numbers: 2, 5, 7, 71, 571, 857, 2857, 28571.

\bigskip

\noindent But unlike cyclic primes, the set of subcyclic primes is always
limited.

\bigskip

\noindent Cyclic and subcyclic primes can be formed from any prime number P
in a given numeric system N. This is a consequence of the fact that every prime number can be at least full repntend prime in some numeric system.

\paragraph{7. Connections between cyclic primes formed from one prime number in different numeric systems} ~\\

\noindent As stated earlier, a cyclic prime can be derived from any prime
number in some numeric systems.

\bigskip

\noindent For each particular cyclic prime, there is a prime P and a numeric system N, with which a cyclic prime was formed.

\bigskip

\noindent However, cyclic primes formed from one P, but in different numeric systems, can exhibit the properties of interconnection.

\bigskip

\noindent For example, in decimal notation, we get prime numbers from a cyclic number 142857.
And in the 40th numeric system, we get prime numbers from the cyclic number 5SMYBH (which corresponds to the sequence of numbers 5, 28, 22, 34, 11, 17).

\bigskip

\noindent However, if we take a prime number that originally looks like H5SMYBH at the 40 numeric system, and translate it into the decimal numeric system, we will see some pattern: 70217142857.

\bigskip

\noindent The least significant bits will correspond to the formation of cyclic primes, but in the most significant digits will be deviations.

\bigskip

\noindent The same behavior is inherent in all cyclic primes found in the 40th system. 

\bigskip

\noindent And in general, a similar pattern will persist for all primes found in certain numeric systems.

\noindent Here is an example of prime numbers found in decimal numeric system from prime P=7:

\bigskip

1) 1428571

2) 71428571

3) 7142857142857

4) 571428571428571

5) 1428571428571428571428571

6) 28571428571428571428571428571

7) 7142857142857142857142857142857

8) 2857142857142857142857142857142857

9) 42857142857142857142857142857142857142857

\bigskip

\noindent Same numbers in the 40th numeric system:

\bigskip

1) MCYB

2) Ra2YB

3) 13NYIMYBH

4) 277Sb5SMYB

5) 1D8TJS2CYBH5SMYB

6) GP98QAT0SMYBH5SMYB

7) 2NbRO471EIMYBH5SMYBH

8) PdGa11UDOPSMYBH5SMYBH

9) 3WAEQ3OR61AQVH5SMYBH5SMYBH

\bigskip

\noindent Here is an example of cyclic prime for P=7 from the 40th numeric system:
\bigskip

1) H5SMYBH

\bigskip

2) Next number is big enough - 77 digits long, it starts from 5SMYBH, and then repeats till 'B':

\noindent {\footnotesize 5SMYBH5SMYBH5SMYBH5SMYBH5SMYBH5SMYBH5SMYBH5SMYBH5SMYBH5SMYBH5SMYBH5SMYBH5SMYB}

\bigskip

\noindent And these are the same numbers but represented in decimal. First one contains 11 decimal digits.

\bigskip

1) 70217142857

\bigskip

\noindent Next one is big enough so its splited into parts, its cyclic part is 12 cycles long, and the whole number is 123 digits long:

\bigskip

2) 3262280440470765442418939358741703168874849426...

...28571428571428571428571428571428571428571428571428571428571428571428571428571

\bigskip

\noindent The rest of cyclic primes from the 40th numeric system are longer, and they all have the same cyclic part that could be seen in decimal.

\bigskip

\noindent If we consider the manifestation of this property in the decimal numeric system for a simple
number P = 7, then the numeric systems in which we will observe such properties can be described
formulas below.

\bigskip

\noindent For the base numeric system N, then to find the numeric systems in which we find the described property, we need to add alternately N * 3, then N * 4.

\bigskip

\noindent We see this behavior in the 40th, 80th, 110th, 150th numeric systems.

\bigskip

\noindent \textbf{Formula for calculating numeric systems forming links with the original system calculus for a given prime}

\bigskip

\noindent Above, we considered the prime number P = 7 and the numeric system N = 10. Generalized the formula for finding related numeric systems can be written as follows:

\begin{equation}
Ns(i) = N + 3 * N * i + ((i + 1) mod 2) * i * N * 4
\end{equation}

\bigskip

\noindent Where i is a non-negative integer. And i = 0 matches the first system where the prime is a full reptend prime.

\bigskip

\noindent However, this formula is only suitable for decimal notation.

\bigskip

\noindent If we try to write out similar formulas for at least several other numeric systems, it turns out that they are not so easy to generalize.

\bigskip

\noindent Formula for N = 3, 10, 17, 31, 38, 59: 

\begin{equation}
Ns (i) = N + 3 * N * i + ((i + 1) mod 2) * i * N
\end{equation}

\bigskip

\noindent Formula for N = 5, 19, 26, 33, 47, 61: 

\bigskip

\begin{equation}
Ns (i) = N + N * i + ((i + 1) mod 2) * i * 5 * N
\end{equation}

\noindent Formula for N = 12: 

\begin{equation}
Ns (i) = N + N * i + ((i + 1) mod 2) * i * 5 * N
\end{equation}

\noindent N = 40 refers to the group formed from N = 10.
The same is true for N = 24, it is also formed from N = 12.

\bigskip

\noindent The difference between the formed groups lies in the fact that related cyclic numbers begin to form in the numeric systems lower than the initial investigated N.

\bigskip

\noindent For example, we are examining the 40th numeric system, and we will meet its patterns in the decimal numeric system. Thus, just as the cyclic primes obtained in the decimal system exhibit a pattern in the 40th, so do the cyclic primes obtained in the 40th numeric system, they show the patterns in the decimal numeric system.
The same is true for the 12th and 24th numeric systems.

\bigskip

\noindent Despite the fact that many numeric systems form the same formulas, others are still different, such as 12.

\bigskip

\noindent So, there are in themselves cyclic primes that are formed in any numeric system, where the original prime is the full reptend. They can be interconnected with cyclic primes in other numeric systems, while some numbers have a relationship only with the numeric systems above, while others also have a connection with the numeric systems below, as in the case with the 40th and decimal numeric systems.

\bigskip

\noindent We can observe a similar construction for P = 5, there are repeated coefficients with respect to the first numeric system in the group.

\bigskip

\noindent And in the case of P = 17, everything becomes much more complicated, you can see that the steps are always equal to base, base * 2,
base * 4, however their alternation seems to be changing.

\bigskip

\noindent Although it turned out to be difficult for me to directly derive formulas right away, it is possible to show that cyclic primes can be formed from each prime number.

\paragraph{List of observations requiring proof or refutation} ~\\

\begin{enumerate}
\item Decomposition of a periodic fraction into an infinite number of geometric series;
\item The presence of cyclic primes in each numeric system, and the fact that they form a potentially
endless set;
\item Connections between cyclic primes in different numeric systems;
\end{enumerate}

\bigskip

\noindent Proof or disproof of any of these hypotheses is highly appreciated! Any information regarding full reptend prime or cyclic numbers is also interesting.

\bigskip

\noindent Most of the calculations and visualizations were done with a little utility. Its code is far from perfect, but it helps to study patterns visually:

\bigskip

\textbf{https://github.com/eversearch/cyclicprime}

\bigskip

\paragraph{Acknowledgments} ~\\

\noindent I wish to show a great gratefulness for the people who helped me to review and edit this work.

\bigskip

1. Eugene Yakshin

2. Nikita Gabdulin

3. Nikolay Artemenkov

4. Oleg Zhukov

\bigskip

\noindent

\paragraph{References} ~\\

\noindent 1. John Leslie "The Philisophy of Arithmetic: Exhibiting a Progressive View on the Theory and Practice of
Calculation", 1820

\bigskip

\noindent 2. Leonard Eugene Dickson "History of the Theory of Numbers", Volume 1, 1952

\bigskip

\noindent 3. David Wells "The Penguin Dictionary of Curious and Interesting Numbers", 1986

\bigskip

\noindent 4. John Horton Conway, Richard Kenneth Guy "The Book of Numbers", 1996

\bigskip

\bigskip

\noindent Kutsenko Konstantin 

\noindent Email: ssidein@gmail.com

\end{document}